\documentstyle[12pt,euscript]{amsart}
\def\mysavedown#1{\edef\mysubs{\mysubs#1}}
\def\mysaveup#1{\edef\mysups{\mysups#1}}
\def\mydown#1{{\mytensor}_{\vphantom{\mysubs}#1}}
\def\myup#1{{\mytensor}^{\vphantom{\mysups}#1}}
\def\tensor#1#2{
  #1
  \def\mytensor{\vphantom{#1}}
  \def\mysubs{\relax}
  \def\mysups{\relax}
  \let\down=\mysavedown
  \let\up=\mysaveup
  #2
  \let\down=\mydown
  \let\up=\myup
  #2
  }
\renewcommand{\phi}{\varphi}
\renewcommand{\epsilon}{\varepsilon}

\newcommand{\Int}{\operatorname{Int}}

\newcommand{\R}{\bold R}
\newcommand{\cross}{\mathbin{\times}}
\newcommand{\dm}{{\partial M}}
\newcommand{\dminfty}{{\partial_\infty M}}
\newcommand{\dmo}{{\partial_0 M}}
\newcommand{\mbar}{{\overline M}}
\newcommand{\gbar}{{\bar g}}
\renewcommand{\hbar}{{\bar h}}
\newcommand{\ghat}{{\widehat g}}
\newcommand{\gtw}{{\widetilde g}}
\newcommand{\gtwsub}{{\tilde g}}
\newtheorem{theorem}{Theorem}[section]
\newtheorem{lemma}[theorem]{Lemma}
\newtheorem{proposition}[theorem]{Proposition}

\newtheorem{bigtheorem}{Theorem}[section]
\theoremstyle{definition}

\numberwithin{equation}{section}
\let\scr=\EuScript

\def\crn#1#2{{\vcenter{\vbox{
        \hbox{\kern#2pt \vrule width.#2pt height#1pt
           }
          \hrule height.#2pt}}}}

\def\Del{\nabla}

\begin{document}
\title
{The spectrum of an asymptotically hyperbolic Einstein manifold}
\author{John M. Lee}
\thanks{Research supported in part by National
Science Foundation grant DMS 91-01832.}
\address{%
Department of Mathematics, GN-50\\
  University of Washington \\
  Seattle, WA 98195 }
\email{lee@@math.washington.edu}
\maketitle

\begin{abstract}
This paper relates the spectrum of the scalar Laplacian of an
asymptotically hyperbolic Einstein metric to the conformal geometry of its
``ideal boundary'' at infinity.  It follows from work of R. Mazzeo that the
essential spectrum of such a metric on an $(n+1)$-dimensional manifold is
the ray $[n^2/4,\infty)$, with no embedded eigenvalues; however, in general
there may be discrete eigenvalues below the continuous spectrum.  The main
result of this paper is that, if the Yamabe invariant of the conformal
structure on the boundary is non-negative, then there are no such
eigenvalues.  This generalizes results of R. Schoen, S.-T. Yau, and D.
Sullivan for the case of hyperbolic manifolds.
\end{abstract}

\section {Introduction}\label{introsection}

One of the most important global invariants of a complete, non-compact
Riemannian manifold $(M,g)$ is the infimum of the $L^2$ spectrum of its
Laplacian $\Delta_g=d^*d$, denoted $\lambda_0(g)$.  In spite of its
importance and the simplicity of its definition, $\lambda_0(g)$ is
notoriously hard to compute except in simple examples.  It reflects in a
subtle way the interaction among the topology of the manifold, the
curvature of the metric, and the asymptotic behavior of curvature and
volume at infinity.

Much more is known about $\lambda_0(g)$ for hyperbolic manifolds
(Riemannian manifolds with constant sectional curvature equal to $-1$).
For example, on $(n+1)$-dimensional hyperbolic space the Laplacian has
purely continuous spectrum, consisting of the ray $[n^2/4,\infty)$, so
$\lambda_0(g) = n^2/4$.  Complete hyperbolic manifolds, which are quotients
of the unit ball $B^{n+1}$ by Kleinian groups, can in general have both
discrete and continuous spectra.  In particular, for a hyperbolic manifold
$(B^{n+1}/\Gamma,g)$ that is {\it geometrically finite} (i.e.\ $\Gamma$ has
a fundamental domain bounded by finitely many geodesic hyperplanes), the
continuous spectrum is still $[n^2/4,\infty)$, but there may be finitely
many discrete eigenvalues below $n^2/4$ \cite{LP}.  In this case, the key
to understanding $\lambda_0(g)$ is the {\it limit set} $\Lambda(\Gamma)$,
which is by definition the set of limit points in the sphere $S^n=\partial
B^{n+1}$ of the action of $\Gamma$ on $\overline{B^{n+1}}$: in
\cite{Sullivan}, D.  Sullivan showed that if $(B^{n+1}/\Gamma,g)$ is
geometrically finite, there are eigenvalues below $n^2/4$ if and only if
$\Lambda(\Gamma)$ has Hausdorff dimension $d>n/2$, and in that case
$\lambda_0(g) = d(n-d)$.

On the other hand, recent work of R. Schoen and S.-T. Yau on conformally
flat manifolds sheds light on $\Lambda(\Gamma)$ from a different direction.
When $M=B^{n+1}/\Gamma$ is geometrically finite and without cusps, the
action of $\Gamma$ on $B^{n+1}$ extends to a free and properly
discontinuous action on the open subset
$\Omega=S^n\setminus\Lambda(\Gamma)$ of the sphere; the quotient
$S=\Omega/\Gamma$ is the ``ideal boundary'' of $M$, and carries a natural
conformal structure $[\ghat]$ which is locally conformally flat.  In
\cite{SY}, Schoen and Yau studied locally conformally flat manifolds with
the aim of understanding the behavior of the {\it Yamabe invariant} $\scr
Y[\ghat]$, which is the infimum of the total scalar curvature functional
$\int_S R_\gtwsub dV_{\gtwsub}$ over unit-volume metrics in the conformal
class $[\ghat]$.  In particular, they were able to characterize those
locally conformally flat structures with non-negative Yamabe invariant:
they showed that a compact, locally conformally flat manifold $(S,[\ghat])$
with $\scr Y[\ghat]\ge 0$ can always be realized as a quotient of a
simply-connected domain $\Omega\subset S^n$ by a Kleinian group; and among
such quotients, the ones with $\scr Y[\ghat]\ge 0$ are exactly those for
which the limit set has Hausdorff dimension $d \le (n-2)/2$.

Combining the implications $\scr Y[\ghat]\ge 0 \implies d \le
(n-2)/2\implies d \le n/2\implies \lambda_0(g) =n^2/4$ yields a direct
relation between the global Riemannian geometry of a hyperbolic manifold
and the global conformal geometry its ideal boundary: if $(M,g)$ is a
complete, non-compact, geometrically finite, $(n+1)$-dimensional hyperbolic
manifold without cusps, and the conformal structure on its ideal boundary
has non-negative Yamabe invariant, then $\lambda_0(g)=n^2/4$.

The purpose of this paper is to extend this result to a much more general
class of Riemannian and conformal structures, where there is no group
action to mediate between the two structures.  Just as many locally
conformally flat manifolds can be realized as ideal
boundaries of complete hyperbolic manifolds, many compact conformal
manifolds can be realized as ``conformal infinities'' of complete,
asymptotically hyperbolic Einstein metrics.  Suppose $\mbar$ is a compact
manifold-with-boundary, and $g$ is a complete Riemannian metric on the
interior $M$ of $\mbar$.  One says that $g$ is {\it conformally compact}
if, for any smooth defining function $\rho$ (i.e.\ a smooth function on
$\mbar$ that is positive in $M$ and vanishes to first order on $\dm$),
$\gbar= \rho^2 g$ extends continuously to a positive definite metric on
$\mbar$; it is {\it conformally compact of order $C^{m,\alpha}$} if
$\gbar\in C^{m,\alpha}(\mbar)$.  The conformal class of the boundary metric
$\ghat = {\left.\gbar\right|}_{T\dm}$ is invariantly defined, and is called
the {\it conformal infinity} of $g$.  For example, if $B^{n+1}/\Gamma$ is a
geometrically finite hyperbolic manifold without cusps, then its hyperbolic
metric is conformally compact, and its conformal infinity is the induced
locally conformally flat structure on the boundary manifold
$\Omega/\Gamma$.

If $(M,g)$ is conformally compact of order at least $C^{2,0}$ and $\rho$ is
any defining function, a straightforward computation as in
\cite{Mazzeo-Hodge} shows that the sectional curvatures of $g$ approach
$-|d\rho|^2_{\gbar}$ at $\dm$.  Accordingly, one says $g$ is {\it
asymptotically hyperbolic of order $C^{m,\alpha}$} if $g$ is conformally
compact of order $C^{m,\alpha}$, $m+\alpha\ge 2$, and $|d\rho|^2_{\gbar}=1$
along $\dm$.

In \cite{GL}, C. R. Graham and I showed that every conformal structure on
$S^n$ which is sufficiently close to that of the round metric is the
conformal infinity of an Einstein metric on $B^{n+1}$ close to the
hyperbolic one.  More recently \cite{fred}, I generalized this result to
conformal structures close to the conformal infinity of any negatively
curved, asymptotically hyperbolic Einstein metric.  It can be shown that
these Einstein metrics are unique within the class of metrics sufficiently
close to the original Einstein metric; they may very well be globally
unique.  Thus at the very least, a large open subset of the conformal
classes on the sphere have concrete realizations as conformal infinities of
Einstein metrics, and to date no one has found any obstructions to a
conformal class on a boundary having such a realization.

It has been shown by R. Mazzeo \cite{Ma3} that the continuous spectrum of
the Laplacian on an asymptotically hyperbolic manifold consists of the ray
$[n^2/4,\infty)$ with no embedded eigenvalues, just as in the hyperbolic
case.  When the asymptotically hyperbolic metric is Einstein, one might
expect a relationship between its spectrum and global conformal invariants
of its conformal infinity.  This expectation is justified by the following
theorem, which is the main result of this paper.

\begin{bigtheorem}\label{maintheorem}
Let $\mbar$ be a compact $(n+1)$-dimensional manifold-with-boundary.
Suppose $g$ is a smooth Einstein metric on $M=\Int \mbar$, which is
asymptotically hyperbolic of order $C^{3,\alpha}$ and has smooth conformal
infinity $[\ghat]$.  If $\scr Y[\widehat g]\ge 0$, then
$\lambda_0(g)=n^2/4$; in particular, $\Delta_g$ has no $L^2$
eigenfunctions.
\end{bigtheorem}

This is proved by constructing a positive function $\phi$ satisfying
$(\Delta_g - n^2/4)\phi\ge 0$, since the existence of a positive
supersolution to $\Delta_g-\lambda$ implies $\lambda\le\lambda_0(g)$
\cite{CY}.  The strategy for finding a suitable test function $\phi$ is
first to prove the existence of positive solutions to the eigenvalue
equation $(\Delta_g+n+1)u=0$ with prescribed growth at infinity, and then
to show that, for an appropriate choice of asymptotic growth, a power of
$u$ satisfies the desired estimate.  The hardest step is a gradient
estimate for $u$ of the form $|du|_g^2\le u^2$.  This estimate proceeds in
two steps.  First we show that the Einstein condition implies that
$|du|^2_g-u^2$ is subharmonic.  Then, by analyzing the asymptotics of $u$
near the boundary, we show that $|du|^2_g-u^2$ has a continuous extension
to the boundary, and, under the assumptions that $g$ is Einstein and $\scr
Y[\ghat]\ge 0$, can be made non-positive there by a judicious choice of
$u$.

It should be mentioned that the technical details of the asymptotic
analysis are made somewhat more complicated by the fact that we are
assuming only rather weak regularity for $\gbar$; if one knew that $\gbar$
were smooth up to $\dm$, one could work in Fermi coordinates for $\gbar$
and the argument would be a little more straightforward.  However, from
asymptotic computations (cf.\ \cite{FG}) we know that smoothness has to
fail generically in even dimensions, since a logarithmic term
appears in the asymptotic expansion for the Einstein equation.  The
Einstein metrics shown to exist in \cite{fred} all have at least
$C^{3,\alpha}$ regularity when the conformal infinity is sufficiently
regular, so Theorem \ref{maintheorem} applies to them.  It is reasonable to
conjecture that any asymptotically hyperbolic Einstein metric with smooth
conformal infinity has (in suitable coordinates) an asymptotic expansion in
powers of $\rho$ and $\log \rho$; since the expansion is regular
up to a term of the form $\rho^{n}\log \rho$ in even dimensions and to all
orders in odd dimensions, the regularity hypothesis of Theorem
\ref{maintheorem} would always be satisfied.

Already in the hyperbolic case one can see that Theorem \ref{maintheorem}
is not sharp: using the results of Sullivan and Schoen/Yau quoted above,
one need only construct an appropriate Kleinian group whose limit set has
dimension $d$ satisfying $(n-2)/2 < d < n/2$ to obtain an example in which
$\scr Y[\ghat] < 0$ but still $\lambda_0(g) = n^2/4$.  Many such examples
are easily obtained by considering compact $d$-dimensional hyperbolic
manifolds $B^d/\Gamma$, and extending $\Gamma$ to act on $B^{n+1}$ in the
obvious way; in this case the limit set is all of the equatorial
$d$-sphere, and when $n$ is odd, $d$ can be chosen strictly between
$(n-2)/2$ and $n/2$.

In these examples, the $L^2$ eigenfunctions seem to arise in some sense
from the topology of $M=B^{n+1}/\Gamma$, reflected in the size of the group
$\Gamma$.  One might be tempted to conjecture that the converse to Theorem
\ref{maintheorem} is true at least for conformally compact Einstein metrics
on the ball.  However, this is also apparently not the case.  H. Pedersen
\cite{Pedersen} constructed a family of asymptotically hyperbolic Einstein
metrics on $B^4$ depending on a real parameter $t$, whose conformal
infinities are the Berger metrics $\ghat_t=\sigma_1^2 +\sigma_2^2 +
t\sigma_3^2$ (written in terms of the standard left-invariant coframe
$\{\sigma_j\}$ for $S^3=SU(2)$).  As $t\to+\infty$, $\scr Y[\ghat_t]$
decreases to $-\infty$, while the Einstein metrics converge uniformly on
compact sets to the complex hyperbolic metric, so $\lambda_0(g_t)\searrow
2$.  Numerical studies performed by my student John Roth suggest that, as
$t$ increases, the eigenvalues below $n^2/4=9/4$ first appear some time
{\it after}\/ $\scr Y[\ghat_t]$ becomes negative, so there is a range of
values of $t$ for which the converse to Theorem \ref{maintheorem} is
violated.

There are various ways in which Theorem \ref{maintheorem} could be
sharpened: for example, by finding a necessary and sufficient condition on
$[\ghat]$ for $\lambda_0(g)$ to be equal to $n^2/4$, or, even better, by
finding a formula for $\lambda_0(g)$ in terms of conformal invariants of
$[\ghat]$.  Similarly, one might hope for a formula for $\scr Y[\ghat]$
in terms of Riemannian invariants of $g$, or at least a necessary and
sufficient condition on $g$ for $\scr Y[\ghat]\ge 0$.  A useful first step
would be to find a reasonable generalization of the Hausdorff dimension of
$\Lambda(\Gamma)$ that makes sense for arbitrary conformal structures, or
at least for those that are conformal infinities of Einstein metrics.

The outline of the paper is as follows: \S\ref{defsection} presents basic
definitions; \S\ref{analytic-section} develops the analytic results about
the Laplacian on asymptotically hyperbolic manifolds that are needed to
solve $(\Delta_g+n+1)u=0$; \S\ref{testfcn-section} shows how to use a
solution $u$ to construct a supersolution to $\Delta_g-n^2/4$, and shows
that $|du|_g^2-u^2$ is subharmonic; and in \S\ref{asymptotic-section} the
asymptotic expansion of $u$ is computed and used to show that
$|du|_g^2-u^2\le 0$ on the boundary.  The Appendix is devoted to a proof of
an elementary lemma regarding H\"older regularity of solutions to ODE's
that I was unable to find in the literature.

I am indebted to Rafe Mazzeo for a number of useful conversations regarding
the ideas in this paper, most particularly for introducing me to the
literature on hyperbolic manifolds and for making a crucial suggestion
regarding the ODE lemma in the Appendix.

\section {Definitions}\label{defsection}

In this section we define our terms, and introduce the metrics and function
spaces we will be be working with.  Throughout this paper, smooth will
always mean $C^\infty$.

Let $\overline M$ be a smooth, compact, $(n+1)$-dimensional
manifold-with-boundary.  For technical reasons that will appear in the
proof of Proposition \ref{laplace-s-prop} below, it is convenient to assume
in general that $\dm$ is the union of two disjoint closed submanifolds,
labeled $\dmo$ and $\dminfty$, with $\dminfty$ non-empty.  We will
write $M =\mbar\setminus \dminfty$, which is a non-compact manifold with
boundary.  We call $\dmo$ the {\it inner boundary} of $M$ and $\dminfty$
its {\it boundary at infinity}.  By {\it background coordinates} for $M$,
we will always mean a smooth coordinate chart for an open subset of
$\mbar$.

A {\it defining function} for $M$ will be a real-valued function $\rho\in
C^1(\mbar)$ which vanishes to first order on $\dminfty$ and is positive
elsewhere.  A {\it smooth defining function} will always mean a defining
function that is smooth on $\mbar$.

As in the introduction, a smooth metric $g$ on $M$ is said to be {\it
conformally compact} if, for some (hence any) smooth defining function
$\rho$, the tensor $\gbar=\rho^2 g$ has a continuous extension to $\mbar$,
and is positive definite there.  We say $g$ is conformally compact of order
$C^{m,\alpha}$ if $\gbar\in C^{m,\alpha}(\mbar)$.  (Note that $g$ is always
smooth up to the inner boundary.)  We say $g$ has {\it smooth conformal
infinity} if the induced boundary metric $\ghat=
{\left.\gbar\right|}_{T\dminfty}$ on $\dminfty$ is smooth for any smooth
defining function $\rho$.

Let $u$ be a function which is $m$-times continuously differentiable on
$M$.  For $s\in\R$ and $0\le\alpha\le 1$, we will define weighted H\"older
norms $\|u\|^{(s)}_{m,\alpha}$ as follows.  First, in the special case in
which $M$ is a smoothly bounded open subset of $\R^{n+1}$, we define the
norms as in \cite{GL}:
\begin{align*}
\|u\|^{(s)}_{m,0} &:= \sum_{l=0}^m\ \sum_{|\gamma|=l} \| d^{-s+l}
\partial^{\gamma} u \| _{L^\infty},\\
\|u\|^{(s)}_{m,\alpha} &:= \|u\|^{(s)}_{m,0} +
\sum_{|\gamma|=m} \
\sup_{x,y}
\left[ \min (d_x^{-s+m+\alpha}, d_y^{-s+m+\alpha})\frac
{|\partial^\gamma u(x) - \partial^\gamma u(y)|}
{|x-y|^\alpha}\right],
\end{align*}
where $d_x$ is the Euclidean distance from $x$ to $\dminfty$, and for a
multi-index $\gamma$, $\partial^\gamma = \partial^{|\gamma|}/\partial
x^\gamma$.  In the more general case of a manifold-with-boundary, the same
norms can be defined using a covering by background coordinate charts and a
subordinate partition of unity in the usual way.

In either case, let $\Lambda^s_{m,\alpha}(M)$ denote the Banach space of
$m$-times continuously differentiable functions on $M$ for which
$\|u\|^{(s)}_{m,\alpha}$ is finite.  We also let $\widetilde\Lambda
^s_{m,\alpha}(M)$ denote the closed subspace of $\Lambda ^s_{m,\alpha}(M)$
consisting of functions that vanish on $\dmo$.  By
\cite[Proposition 3.3(13)]{GL} these spaces are independent of choices of
coordinates or partition of unity.  Near $\dmo$ these are just functions
which are $C^{m,\alpha}$ up to the boundary (and, in the case of the
$\widetilde\Lambda$ spaces, satisfy homogeneous Dirichlet boundary
conditions there), while near $\dminfty$ they satisfy degenerate estimates
that are weaker than the usual H\"older estimates.  The properties of these
spaces are summarized in \cite[Prop.\ 3.3]{GL}; one important property we
will use repeatedly is the following: if $m\ge 1$, for any smooth defining
function $\rho$,
\begin{displaymath}
\Lambda^s_{m,\alpha}(M) = \rho^s\Lambda^0_{m,\alpha}(M) =
\{ \rho^s u: u\in \Lambda^0_{m,\alpha}(M)\}.
\end{displaymath}

Throughout this paper, we will use the summation convention and classical
index notation, with covariant derivatives indicated by indices preceded by
a semicolon.  Barred quantities will denote covariant derivatives and
curvature invariants computed with respect to $\gbar$, and unbarred
quantities with respect to $g$, with the indicated metric and its inverse
used in each case to lower and raise indices.  We will usually omit the
semicolon for covariant derivatives of a scalar function.  For example,
\begin{displaymath}
\bar \rho^{i} = \gbar^{ij}\rho_j = \gbar^{ij}\partial\rho/\partial x^j
\end{displaymath}
are the components of the $\gbar$-gradient of $\rho$.

\section {Analytic Preliminaries}\label{analytic-section}

The purpose of this section is to present some basic analytic facts about
the Laplace operator on conformally compact manifolds.  These results can
be derived, for example, from the parametrix construction for elliptic edge
operators due to Rafe Mazzeo \cite{Mazzeo-edge}; however, the results we
need follow in a straightforward way from the analysis in \cite{GL}, so I
have chosen to present direct proofs.

Throughout this section $(M,g)$ will be asymptotically hyperbolic of order
$C^{k,\beta}$, $k\ge 1$, possibly with inner boundary; $\rho$ will be a
fixed smooth defining function; and $\gbar = \rho^2 g$.  Let
$\Delta_g=d^*d$ denote the scalar Laplacian of $g$.  By a slight
generalization of the arguments in \cite[\S3]{GL}, we obtain the following
result:

\begin{lemma}\label{laplace-0-lemma}
For a non-negative integer $m$ and real number $0<\alpha<1$ with
$m+1+\alpha \le k+\beta$, let $h\in\Lambda^0_{m,\alpha}(M)$ be a
real-valued function such that $h\ge\delta>0$ on $M$; and let $X$ be a
vector field on $M$ whose coefficients in any background coordinate chart
are in $\Lambda^1_{m,\alpha}(M)$.  Then
\begin{displaymath}
\Delta_g + X + h\colon \widetilde\Lambda^0_{m+2,\alpha}(M) \to
                                 \Lambda^0_{m,\alpha}(M)
\end{displaymath}
is an isomorphism.
\end{lemma}

\begin{pf}
In the special case that $X=0$, $\dmo = \emptyset$, and $\gbar$ is smooth,
this follows directly from Propositions 3.4, 3.7, and 3.8 (and the remark
following Proposition 3.8) in \cite{GL}.  The generalization to the present
case is straightforward: First, the presence of the first-order term $X$
has no effect on the argument, as one can easily check; in particular, the
test function estimate $(\Delta_g + X + h)\phi \ge \delta\phi$ still holds
with $\phi=1$, and the arguments of \cite[Prop.\ 3.8]{GL} go through
unchanged.  Second, the presence of the inner boundary $\dmo$ can be
handled by applying the usual boundary regularity theory and Schauder
estimates for uniformly elliptic Dirichlet problems.  Finally, if $\gbar$
is only in $C^{k,\beta}(\mbar)\subset C^{m+1,\alpha}(\mbar)$, the same
argument goes through; the only point that requires attention is to check
that the coefficients of the $r$-th order terms in $\Delta_g$ are in
$\Lambda^r_{m,\alpha}(M)$, which is just what is needed to insure that
there is a uniform constant $C$ in the estimate of \cite[Prop.\ 3.4]{GL}.
\end{pf}

Using this lemma, we will prove a fundamental weighted isomorphism result
for the Laplacian plus a constant.  First we need to examine how the
Laplacian behaves when conjugated by a power of $\rho$.

\begin{lemma}\label{conjugate-laplacian-lemma}
Suppose $0<\alpha<1$ and $m+1+\alpha \le k+\beta$.  For $\epsilon>0$, let
$M_\epsilon =\{ x\in M: \rho(x) \le\epsilon\}$.  For real numbers $\kappa$
and $s$,
\begin{equation}\label{epsilon-iso}
\rho^{-s} (\Delta_g+\kappa)\rho^s\colon
\widetilde\Lambda^0_{m+2,\alpha}(M_\epsilon)
\to \Lambda^0_{m,\alpha}(M_\epsilon)
\end{equation}
is an isomorphism whenever
\begin{equation}\label{s-condition}
\left|s - \frac n 2\right|^2 < \frac {n^2} 4 + \kappa
\end{equation}
and $\epsilon$ is sufficiently small.
\end{lemma}

\begin{pf}
Computing as in \S2 of \cite{GL}, we obtain
\begin{align*}
\rho_{ij} &= \partial_i \partial_j\rho -
\Gamma_{ij}^k\partial _k\rho\\
&= \partial_i\partial_j \rho - (\overline\Gamma_{ij}^k -\rho^{-1}
(\delta_i^k\rho_{j} +
\delta_j^k\rho_{i} - \gbar_{ij}\bar\rho^{k}))\rho_k\\
&= \bar\rho_{ij} + \rho^{-1}(2 \rho_i\rho_j - \rho_k\bar\rho^k
\gbar_{ij}),
\end{align*}
and so
\begin{displaymath}
\rho_i{}^i = \rho^2 \gbar^{ij} \rho_{ij} =
\rho^2 \bar\rho_{i}{}^{i}-(n-1) \rho \rho_i\bar\rho^i.
\end{displaymath}
Using this, we obtain
\begin{equation}\label{conjugate-laplacian}
\begin{aligned}
\rho^{-s} (\Delta_g+\kappa)(\rho^s u) &= \rho^{-s} ( -(\rho^s u)_{;i}{}^i +
\kappa \rho^s u)\\
&= -u_i{}^i - s \rho^{-1} \rho_i{}^i u - 2 s \rho^{-1} \rho^i u_i - s (s-1)
\rho^{-2} \rho_i\rho^i u + \kappa u\\
&= -u_i{}^i - s\rho  \bar\rho_{i}{}^{i} u + s(n-1)
\rho_i\bar\rho^i u - 2 s \rho \bar\rho^i u_i - s(s-1)
\rho_i\bar\rho^i u +\kappa u \\
&= (\Delta_g + X + h) u,
\end{aligned}
\end{equation}
where
\begin{align*}
X u &= - 2 s \rho \bar \rho^i u_i,\\
h &=  \kappa + s(n-s)|d\rho|_{\gbar}^2 +
s\rho  \Delta_{\gbar}
\rho.
\end{align*}

Since $|d\rho|^2_{\gbar}\to 1$ and $\rho\Delta_{\gbar}\rho\to 0$ at
$\dminfty$, it follows that $h$ is uniformly positive near $\dminfty$
provided $\kappa+s(n-s)> 0$, which is equivalent to (\ref{s-condition}).
Choosing $\epsilon$ sufficiently small, we see that $\rho^{-s}
(\Delta_g+\kappa)\rho^s= \Delta_g+X+h$ satisfies the hypotheses of Lemma
\ref{laplace-0-lemma} on $M_\epsilon$, and the lemma follows.
\end{pf}

\begin{proposition}\label{laplace-s-prop}
Let $\kappa$ be a positive real number, $0<\alpha<1$, and $m+1+\alpha \le
k+\beta$.  Then
\begin{equation}\label{delta+kappa-iso}
\Delta_g + \kappa\colon \widetilde\Lambda^s_{m+2,\alpha}(M)\to
\Lambda^s_{m,\alpha}(M)
\end{equation}
is an isomorphism whenever {\rom(}\ref{s-condition}{\rom)} holds.
\end{proposition}

\begin{pf}
Given $f\in \Lambda^s_{m,\alpha}(M)$, we need to show that there exists a
unique $u\in\widetilde\Lambda^s_{m+2,\alpha}(M)$ such that
\begin{equation}\label{delta+kappa}
(\Delta_g+\kappa) u= f.
\end{equation}
By the closed graph theorem, it follows then that $(\Delta_g+\kappa)$
has a bounded inverse.

The case $s=0$ is just a special case of Lemma \ref{laplace-0-lemma}.
Consider next the case $s<0$.  We can write $f = \rho^s b$ with $b\in
\Lambda^0_{m,\alpha}(M)$.
Then choosing $\epsilon$ small enough that (\ref{epsilon-iso}) is an
isomorphism, we can find $v\in\widetilde\Lambda^0_{m+2,\alpha}(M)
(M_\epsilon)$ satisfying $\rho^{-s} (\Delta_g+\kappa)(\rho^s v)=b$ on
$M_\epsilon$.

Let $\phi\in C^\infty(\mbar)$ be a function that is supported in
$M_\epsilon$ and identically equal to 1 in a neighborhood of $\dminfty$.
Then $\rho^s b-(\Delta_g+\kappa) (\phi\rho^s v)\in \Lambda^0_{m,\alpha}(M)$
since it vanishes near $\dminfty$; therefore by Lemma
\ref{laplace-0-lemma} again there exists
$w\in\widetilde\Lambda^0_{m+2,\alpha}(M)$ satisfying
\begin{displaymath}
(\Delta_g+\kappa)w =   \rho^s b-(\Delta_g+\kappa)(\phi\rho^s v).
\end{displaymath}
It follows immediately that $u= w+\phi\rho^s v$ solves (\ref{delta+kappa}),
and since $\widetilde\Lambda^0_{m+2,\alpha}(M)\subset \widetilde \Lambda
^s_{m+2,\alpha}(M)$, we have $u\in \widetilde \Lambda ^s_{m+2,\alpha}(M)$.

For uniqueness, suppose that $(\Delta_g+\kappa)u=0$ for some $u\in
\widetilde\Lambda^s_{m+2,\alpha}(M)$.  With $\phi$ as above, it is easy to
see that $h=(\Delta_g+\kappa)(\phi u)\in
\Lambda^0_{m,\alpha}(M_\epsilon)$ because it vanishes near $\dminfty$.
Therefore, since (\ref{epsilon-iso}) is an isomorphism for $s=0$ and
$\epsilon$ small, there is a unique $v\in
\widetilde\Lambda^0_{m+2,\alpha}(M_\epsilon)$ satisfying
$(\Delta_g+\kappa)v=h$.  But
$(\Delta_g+\kappa)$ is also injective on
$\widetilde\Lambda^s_{m+2,\alpha}(M_\epsilon)$, so we must have $\phi u = v
\in \widetilde\Lambda^0_{m+2,\alpha}(M_\epsilon)$, so $u$
itself is in $\widetilde\Lambda^0_{m+2,\alpha}(M)$, and by injectivity of
(\ref{delta+kappa-iso}) in
the case $s=0$ we conclude finally that $u=0$.

Now suppose $s>0$, and $f = \rho^s b\in \Lambda^s_{m,\alpha}(M)$.  Since
$\Lambda^s_{m,\alpha}(M)\subset \Lambda^0_{m,\alpha}(M)$, by the $s=0$ case
there exists $u\in \widetilde\Lambda^0_{m+2,\alpha}(M)$ satisfying
(\ref{delta+kappa}).  We need only show that $u\in\widetilde
\Lambda^s_{m+2,\alpha}(M)$, for then injectivity on
$\widetilde\Lambda^0_{m+2,\alpha}(M)$ automatically implies injectivity on
the smaller space $\widetilde \Lambda^s_{m+2,\alpha}(M)$.

Choosing $\epsilon$ small enough that (\ref{epsilon-iso}) is an
isomorphism, and letting  $\phi$ be a cutoff function as before, observe
that $\rho^{-s}(\Delta_g+\kappa)(\phi u)\in
\widetilde\Lambda^0_{m,\alpha}(M_\epsilon)$ since it agrees with
$b=\rho^{-s}f$ near $\dminfty$.  Therefore by Lemma
\ref{conjugate-laplacian-lemma} there exists $v\in
\widetilde\Lambda^0_{m+2,\alpha}(M_\epsilon)$ satisfying
\begin{displaymath}
\rho^{-s} (\Delta_g+\kappa) \rho^s v = \rho^{-s} (\Delta_g+\kappa)(\phi u)
\end{displaymath}
on $M_\epsilon$.  By injectivity of (\ref{delta+kappa-iso}) when $s=0$ it
follows that $\phi u = \rho^s v$, which completes the proof.
\end{pf}

\section{The test function estimate}\label{testfcn-section}

Armed with the isomorphism result of Proposition \ref{laplace-s-prop}, we
will show in this section how to use a power of a solution to the
eigenvalue equation $(\Delta_g+n+1)u=0$ to construct a test function $\phi$
satisfying $(\Delta_g-n^2/4)\phi \ge 0$.  Throughout this section, we will
assume only that $(M,g)$ is asymptotically hyperbolic of order
$C^{1,\alpha}$, $0<\alpha<1$, without inner boundary.

\begin{proposition}\label{uprop}
Let $\rho$ be any smooth defining function.  There exists a unique smooth,
strictly positive function $u$ on $M$ satisfying
\begin{gather}
(\Delta_g+n+1)u = 0\label{laplace+n+1},\\
u-\rho ^{-1} \text{ is bounded.}\label{u-rhoinv}
\end{gather}
\end{proposition}

\begin{pf}
Let $\gbar=\rho^2 g$.  From (\ref{conjugate-laplacian}) with $u=1$, we
obtain
\begin{equation}\label{laplacian-rho-1}
(\Delta_g+n+1) \rho^{-1} = \bar\rho_{i}{}^{i} + (n+1)\rho^{-1}
(1-\rho_i\bar\rho^i).
\end{equation}
The assumption that $\gbar\in C^{1,\alpha}(\mbar)$ guarantees that
$\bar\rho_{i}{}^{i}\in C^{0,\alpha}(\mbar)$.  Moreoever,
$1-\rho_i\bar\rho^i\in C^{1,\alpha}(\mbar)$ and vanishes on $\dminfty$ by
the assumption of asymptotic hyperbolicity, so $\rho^{-1}
(1-\rho_i\bar\rho^i)\in C^{0,\alpha}(\mbar)$.  Therefore,
$(\Delta_g+n+1)\rho^{-1}\in C^{0,\alpha}(\mbar)\subset
\Lambda^0_{0,\alpha}(M)$.  It follows immediately from Proposition
\ref{laplace-s-prop} that there exists a unique function $v\in
\Lambda^0_{2,\alpha}(M)$ such that $u = \rho^{-1} + v$ satisfies
(\ref{laplace+n+1}) and (\ref{u-rhoinv}).

Now, since $g$ is smooth in $M$, $u$ is smooth by local elliptic
regularity.  The strong maximum principle applied to (\ref{laplace+n+1})
shows that $u$ cannot attain a non-positive interior minimum.  Thus, since
$u\to+\infty$ at $\dminfty$, $u$ is strictly positive on $M$.

To prove uniqueness, let $v$ be the difference between two solutions; then
$v$ is a bounded solution to $(\Delta_g+n+1)v=0$.  Any bounded function is
in $\Lambda^{0}_{0,0}(M)$ by definition, and then it follows from
\cite[Proposition 3.4]{GL} that $v\in \Lambda^0_{2,\alpha}(M)$.  Using the
injectivity assertion of Proposition \ref{laplace-s-prop}, we conclude that
$v=0$.
\end{pf}

We will take as our test function $\phi = u^{-s}$, for a real power $s$ to
be determined shortly.  A straightforward computation gives
\begin{equation}\label{phi-estimate}
\frac{\Delta_g\phi}{\phi} = -s\frac {\Delta_g u}{u} - s(s+1) \frac
{|du|_g^2} {u^2} = s(n+1) - s(s+1) \frac {|du|_g^2} {u^2}.
\end{equation}
Since $u$ looks asymptotically like $\rho^{-1}$ near $\dminfty$, it is
reasonable to expect that, asymptotically,
\begin{displaymath}
\frac{|du|_g^2}{u^2} \sim \frac {|\rho^{-2}d\rho|_g^2}{\rho^{-2}} =
\rho^{-2} |d\rho|_g^2 = |d\rho|_{\gbar}^2 \sim 1.
\end{displaymath}
Therefore, near the boundary, $\Delta_g\phi/\phi$ approaches $s(n-s)$, so
the best global estimate we can hope for is that $\Delta_g\phi/\phi \ge
s(n-s)$.  For $s>0$, (\ref{phi-estimate}) shows that this is true precisely
when $|du|_g^2 /u^2 \le 1$ on $M$.  In that case $\phi$ is a supersolution
to $\Delta_g-n^2/4$, so $s(n-s)$ is a lower bound for $\lambda_0(g)$.  This
bound is optimal when $s=n/2$, in which case we get the lower bound $n^2/4$
claimed in Theorem
\ref{maintheorem}.

This argument reduces the proof of Theorem \ref{maintheorem} to the global
gradient estimate
\begin{equation}\label{gradient-estimate}
|du|_g^2 \le u^2
\end{equation}
for the eigenfunction $u$ of Proposition \ref{uprop}.  We will accomplish
this in two steps.  First, in the next proposition, we will show that
$|du|_g^2 - u^2$ is subharmonic when $g$ is Einstein, so it cannot take an
interior maximum; then in the next section we will show that, under the
hypotheses of Theorem \ref{maintheorem}, $|du|_g^2 - u^2 \le 0$ on
$\dminfty$, so that (\ref{gradient-estimate}) holds globally on $M$.

\begin{proposition}\label{subharmonic-prop}
Suppose $g$ is any Riemannian metric with $Rc_g \ge -n g$.  Let $u$ be any
solution to $(\Delta_g+n+1) u= 0$.  Then $|du|_g^2-u^2$ is subharmonic.
\end{proposition}

\begin{pf}
A standard calculation using the Ricci identity and  equation
(\ref{laplace+n+1})
shows that $u_{ij}{}^j = (n+1)u_i + R_{ij}u^j$, and so
\begin{equation}\label{delta-du2-u2}
\begin{aligned}
\Delta_g ( |du|_g^2 - u^2)
&= -(u_i u^i - u^2)_{;j}{}^j\\
&= - 2 u_{ij}{}^j u^i - 2 u_{ij}u^{ij} + 2 u_j{}^j u + 2 u_j u^j\\
&= - 2n u_j u^j - 2 R_{ij} u^i u^j - 2 u_{ij} u^{ij} + 2 (n+1)u^2.
\end{aligned}
\end{equation}
Let $b = \Del_g^2 u + \frac {1}{n+1}(\Delta_g u) g$ denote the traceless
Hessian of
$u$.  Then
\begin{displaymath}
|b|_g^2 = (u_{ij} - u g_{ij})(u^{ij} - u g^{ij}) = u_{ij} u^{ij} -
(n+1)u^2.
\end{displaymath}
Inserting this into (\ref{delta-du2-u2}), we get
\begin{displaymath}
\Delta_g ( |du|_g^2 - u^2) = - 2 \left< (Rc_g+n) du,
du\right>_g - 2 |b|_g^2  \le 0,
\end{displaymath}
which completes the proof.
\end{pf}

\section{Boundary Asymptotics}\label{asymptotic-section}

In this section, we will compute the boundary asymptotics of the
eigenfunction $u$ of Proposition \ref{uprop} on an asymptotically
hyperbolic Einstein manifold.  The main result is Proposition
\ref{boundary-prop}, which gives the boundary values of $|du|_g^2-u^2$. At
the end of the section, we prove Theorem \ref{maintheorem}.

We assume in this section that $(M,g)$ is asymptotically hyperbolic of
order $C^{3,\alpha}$, without inner boundary, and $g$ is Einstein.

\begin{lemma}\label{r-lemma}
Let $\rho$ be any smooth defining function for $M$.  Then there exists
a defining function $r\in C^{2,\alpha}(\mbar)$ such that
\begin{enumerate}
\item\label{rho+rho^2}
$r = \rho + O(\rho^2)$;
\item\label{c2alpha}
$\gtw = r^2 g$ has a $C^{2,\alpha}$ extension to $\mbar$; and
\item\label{|dr|^2=1}
$|dr|_{\gtwsub}^2 \equiv 1$ in a neighborhood of $\dminfty$.
\end{enumerate}
\end{lemma}

\begin{pf}
Writing $\gbar=\rho^2 g$ and $r = e^v \rho$, and computing as in the proof
of Lemma 5.2 of \cite{GL}, we see that (\ref{|dr|^2=1}) is equivalent to
\begin{equation}\label{F-eqn}
2\left<d\rho,dv\right>_{\gbar} + \rho|dv|_{\gbar}^2 = \frac
{1-|d\rho|^2_{\gbar}} {\rho}.
\end{equation}
If $v\in C^{2,\alpha}(\mbar)$ is a solution to (\ref{F-eqn}), then
(\ref{c2alpha}) and (\ref{|dr|^2=1}) will be satisfied, and
(\ref{rho+rho^2}) will follow if in addition $v=0$ on $\dminfty$.

Now (\ref{F-eqn}) is a non-characteristic first-order PDE for $v$, and can
be solved easily in a neighborhood of $\dminfty$ by Hamilton-Jacobi theory.
However, since the standard treatments of Hamilton-Jacobi theory do not
give H\"older regularity of the solution when the equation has H\"older
coefficients, we will go through the construction in some detail.

Let $F\colon T^*\mbar\to \R$ denote the function
\begin{displaymath}
F(x,\xi) = 2\left< d\rho(x),\xi\right>_{\gbar} + \rho(x)|\xi|_{\gbar}^2 -
\frac {1-|d\rho(x)|^2_{\gbar}} {\rho(x)}
\end{displaymath}
for $x\in \mbar$, $\xi\in T_{x}^*\mbar$.  Solving (\ref{F-eqn}) is
equivalent to finding a function $u$ on $\mbar$ such that $F(x,du(x))\equiv
0$.

Since $\gbar\in C^{3,\alpha}(\mbar)$ by hypothesis, the first two terms in
$F$ are $C^{3,\alpha}$ functions on $T^*\mbar$.  Since
$1-|d\rho|_{\gbar}^2$ is in $C^{3,\alpha}(\mbar)$ and vanishes on
$\dminfty$, it follows that $(1-|d\rho|_{\gbar}^2)/\rho\in
C^{2,\alpha}(\mbar)$, and hence $F\in C^{2,\alpha}(T^*\mbar)$.  Thus the
Hamiltonian vector field $X_F$ of $F$ is in $C^{1,\alpha}(T^*\mbar)$, and
by Lemma \ref{ode-lemma} in the Appendix, its flow $\phi\colon{\R}
\cross T^*\mbar \to T^*\mbar$ is a $C^{1,\alpha}$ map wherever it is
defined.

Write $b = (1-|d\rho|_\gbar^2)/\rho$, and let $\psi\colon
\overline{\R^+}\cross\dminfty\to T^*\mbar$ be the $C^{1,\alpha}$
map
\begin{displaymath}
\psi(t,x) = \phi(t, (x,\tfrac12 b(x)d\rho(x))).
\end{displaymath}
In other words, $\psi$ is the flow-out by $X_F$ from the set $\{(x,\tfrac12
b(x)d\rho(x))\in T^*\mbar: x\in \dminfty\}$.  Since $F(x,\tfrac12
b(x)d\rho(x))=0$ for $x\in \dminfty$, a standard argument in symplectic
geometry shows that the image of $\psi$ is a $C^{1,\alpha}$ Lagrangian
submanifold-with-boundary of $T^*\mbar$ contained in $F^{-1}(0)$.  It is
easy to check that it is transverse to the fibers of $T^*\mbar$ along
$\dminfty$, and so near $\dminfty$ it is the image of a closed $1$-form
$\omega\in C^{1,\alpha}(\mbar,T^*\mbar)$.

At least locally, $\omega=dv$ for some function $v\in C^{2,\alpha}(\mbar)$.
Since $\omega = \tfrac12 b d\rho$ along $\dminfty$, any such
function $v$ is locally constant on $\dminfty$; if we require $v=0$ on
$\dminfty$,
then $v$ is uniquely determined, and thus globally defined in a
neighborhood of $\dminfty$.  Extending $v$ arbitrarily to all of $\mbar$ and
setting $r=e^v \rho$ completes the proof of the lemma.
\end{pf}

\begin{lemma}\label{u-asymptotics-lemma}
With $u$ as in Proposition \ref{uprop} and $r$ as in the preceding lemma,
\begin{displaymath}
u = r^{-1} + r v,
\end{displaymath}
where $v$ is a $C^{2,\alpha}$ function on $M$ such that $v$ and $|dv|_g$
have continuous extensions to $\mbar$, with
\begin{displaymath}
v= \frac 1 {4n(n-1)} \widehat R \quad \text{and}\quad |dv|_g= 0 \quad
\text{on $\dminfty$},
\end{displaymath}
where $\widehat R$ is the scalar curvature of $\ghat = {\left.r^2
g\right|}_{T\dminfty}$.
\end{lemma}

\begin{pf}
Writing $\gtw = r^2 g$, the formula for the transformation of the Ricci
tensor under a conformal change of metric (cf., for example,
\cite[p.\ 59]{Besse}) can be written
\begin{displaymath}
R_{ij} = \widetilde R_{ij} + (n-1) r^{-1} \tilde r_{ij} + r^{-1}\tilde
r_k{}^k\gtw_{ij} - n r^{-2} r_k \tilde r^k \gtw_{ij},
\end{displaymath}
where quantities with tildes are computed with respect to $\gtw$.
Using $r_k \tilde r^k = 1$ and multiplying by $r$, near $\dminfty$ the
Einstein equation $R_{ij} =-ng_{ij}$ becomes
\begin{equation}\label{Einstein}
r\widetilde R_{ij} + (n-1) \tilde r_{ij} + \tilde r_k{}^k
\gtw_{ij} = 0.
\end{equation}
We will need two simple consequences of this formula.  First, taking the
trace with respect to $\gtw$, we get
\begin{equation}\label{traced-Einstein}
\tilde r_i{}^i = - \frac 1 {2n} r \widetilde R.
\end{equation}
Second, contracting (\ref{Einstein}) with $\tilde r^i \tilde r^j$ and using
the fact that $\tilde r_{ij} \tilde r^i \tilde r^j = \tfrac12 (r_i \tilde
r^i)_{;j}
\tilde r^j = 0$,
\begin{equation}\label{normal-Einstein}
\tilde r_i{}^i = - r\widetilde R_{ij}\tilde r^i \tilde r^j.
\end{equation}
These last two equations together imply
\begin{equation} \label{R=Rijrirj}
\widetilde R_{ij}\tilde r^i \tilde r^j = \frac 1 {2n} \widetilde R.
\end{equation}

Next we will relate $\widetilde R$ with the scalar curvature $\widehat R$
of the boundary metric $\ghat$.  From (\ref{Einstein}) and
(\ref{normal-Einstein}) it follows that $\tilde r_{ij}=0$ on $\dminfty$, so
$\dminfty = \{r=0\}$ is totally geodesic in the $\gtw$ metric.  Therefore,
at any point $x\in \dminfty$, in terms of an orthonormal basis $\{X_0,
X_\alpha,\ \alpha = 1,\dots,n\}$ for $T_x \mbar$ with $X_0$ normal to
$\dminfty$ and $X_\alpha$ tangent to $\dminfty$, the curvature tensors of
$\gtw$ and $\ghat$ satisfy $\widetilde R_{\alpha\beta\gamma\delta} =
\widehat R_{\alpha\beta\gamma\delta}$.  Consequently the respective scalar
curvatures satisfy
\begin{equation}\label{rtw=rhat+2r00}
\widetilde R = \widetilde R^{\alpha\beta}{}_{\alpha\beta}
+ \widetilde R^{\alpha 0}{}_{\alpha 0}
+ \widetilde R^{0\beta}{}_{0\beta}
= \widehat R + 2 \widetilde R_{00}.
\end{equation}
Since $|dr|_{\gtwsub} = 1$, $\widetilde R_{00} =\widetilde R_{ij}\tilde r^i
\tilde r^j$ on $\dminfty$, and so from (\ref{rtw=rhat+2r00}) and
(\ref{R=Rijrirj}) we conclude that
\begin{equation}\label{rtilde=rhat}
\widehat R = \widetilde R - 2 \widetilde R_{ij}\tilde r^i \tilde r^j =
\frac {n-1} n \widetilde R \qquad \text{on $\dminfty$}.
\end{equation}

 From (\ref{laplacian-rho-1}), near $\dminfty$ we have
\begin{displaymath}
(\Delta_g+n+1) r^{-1} = \tilde r_i{}^i = - \frac 1 {2n} r \widetilde R.
\end{displaymath}
Let $f\in C^\infty(\mbar)$.  Then a computation similar to
(\ref{conjugate-laplacian}) yields
\begin{displaymath}
(\Delta_g+n+1) (r f) = 2nrf + r^2 \big((n-3) \tilde r^i f_i - f \tilde r_i{}^i
- r \tilde f_i{}^i\big).
\end{displaymath}
 From the last three equations, we conclude that, if $f$ is chosen so that
${\left. f \right|}_{\dminfty}=\widehat R$,
\begin{displaymath}
(\Delta_g+n+1) \left(r^{-1} + \frac 1 {4n(n-1)} r f\right) = r h,
\end{displaymath}
for some function $h\in C^{0,\alpha}(\mbar)$ which vanishes on $\dminfty$.
If $0<\epsilon\le\alpha/2$, it can be verified easily that $rh\in
\Lambda^{1+\epsilon}_{0,\epsilon}(M)$, and so from Proposition
\ref{laplace-s-prop} we conclude there is a unique function
$w\in\Lambda^{0}_{2,\epsilon}(M)$ such that $r^{1+\epsilon} w\in
\Lambda^{1+\epsilon}_{2,\epsilon}(M)$ satisfies
\begin{displaymath}
(\Delta_g+n+1)\biggl(r^{-1} + \frac 1 {4n(n-1)} r f +
r^{1+\epsilon}w\biggr) = 0.
\end{displaymath}
Since $r^{-1} - \rho^{-1}$ is bounded, by the uniqueness part of
Proposition \ref{uprop}, we conclude that
\begin{displaymath}
u = r^{-1} +  \frac 1 {4n(n-1)} r f + r^{1+\epsilon}w.
\end{displaymath}
If we set $v = f/(4n(n-1)) + r^\epsilon w$, it follows that
\begin{displaymath}
{\left.v\right|}_{\dminfty} = \frac 1 {4n(n-1)} {\left.f\right|}_{\dminfty} =
\frac 1 {4n(n-1)} \widehat R
\end{displaymath}
and
\begin{displaymath}
|dv|_g = r|dv|_{\gtwsub}\le C( r |df|_{\gtwsub} + r^{1+\epsilon} |dw|_{\gtwsub}
+ |w|
r^{\epsilon} |dr|_{\gtwsub}) \le C' r^{\epsilon}.
\end{displaymath}
This completes the proof.
\end{pf}

\begin{proposition}\label{boundary-prop}
Suppose $\rho$ is any smooth defining function for $M$, and let $\ghat =
{\left.\rho^2 g\right|}_{T\dminfty}$.  If $u$ is the solution to
$(\Delta_g+n+1)u=0$ given by Proposition
\ref{uprop},  then $|du|^2_g-u^2$ has a continuous extension to $\mbar$,
and is equal to $-\tfrac{1}{n(n-1)}\widehat R$ along $\dminfty$.
\end{proposition}

\begin{pf}
Using the asymptotic expansion of $u$ given by the preceding lemma together
with the fact that $|dr|^2_g = r^2$ near $\dminfty$, we
compute
\begin{align*}
|du|_g^2-u^2 &= |(-r^{-2})dr + v dr + r dv|_g^2 - (r^{-1} + r v)^2\\
&=
r^{-4} |dr|_g^2 + v^2 |dr|_g^2 + r^2 |dv|_g^2 - 2 r^{-2} v|dr|_g^2 - 2
r^{-1} \left<dr,dv\right>_g\\
&\qquad + 2 r v \left <dr,dv\right>_g - r^{-2} - 2 v -
r^2 v^2\\
&= - 4 v - 2 r^{-1} \left<dr,dv\right>_g +O(r^2).
\end{align*}
Since $|dv|_g\to 0$ at $\dminfty$, this expression approaches $-4v= -\tfrac
1 {n(n-1)} \widehat R$ at $\dminfty$, which completes the proof.
\end{pf}

\begin{pf*}{Proof of Theorem \ref{maintheorem}}
If $g$ is as in the statement of the theorem, $\rho$ is any smooth defining
function for $M$, and $u$ is as in Proposition \ref{uprop}, Proposition
\ref{boundary-prop} shows that $|du|_g^2 - u^2$ approaches a negative
constant times $\widehat R$ at $\dminfty$.  By the solution of the Yamabe
problem \cite{Schoen}, $\rho$ can be chosen so that $\widehat R$ is
constant on $\dminfty$, and the assumption that $\scr Y[\ghat]\ge 0$
guarantees that this constant is non-negative.  Since $|du|_g^2-u^2$ is
subharmonic by Proposition \ref{subharmonic-prop} and non-positive on
$\dminfty$, it follows from the maximum principle that $|du|_g^2\le u^2$ on
all of $M$, and thus from (\ref{phi-estimate}) the positive function $\phi
= u^{-n/2}$ satisfies $\Delta_g \phi \ge (n^2/4)\phi$ on $M$, so the
spectrum of $\Delta_g$ is bounded below by $n^2/4$.
\end{pf*}

\section*{Appendix: H\"older regularity for ODE's}
\renewcommand{\thetheorem}{A.\arabic{theorem}}
\renewcommand{\theequation}{A.\arabic{equation}}
\setcounter{equation}{0}
\setcounter{theorem}{0}

The purpose of this appendix is to extend the standard regularity theory
for ordinary differential equations to show that the flow of a vector field
with H\"older coefficients has the same H\"older regularity as the
coefficients.  This is probably well-known to PDE experts, but I was unable
to find it in the literature, so I have included a proof of the
$C^{1,\alpha}$ case, which is all that is needed for this paper.

\begin{lemma}\label{ode-lemma}
Let $V$ be a vector field on an open subset in $\R^m$, whose coefficient
functions are of class $C^{1,\alpha}$, $0 < \alpha < 1$.  The solution
$y(x,t)$ to the initial value problem
\begin{align}
\frac{\partial y}{\partial t}(x,t) &= V(y(x,t)),\label{ode}\\
y(x,0) &= x\label{iv}
\end{align}
is a $C^{1,\alpha}$ function of $(x,t)$ where it is defined.
\end{lemma}

\begin{pf}
By classical ODE regularity theory, $y$ is a $C^1$ function of $(x,t)$.
Since $\partial y^i/\partial t = V^i\circ y \in C^{1}\subset C^{0,\alpha}$,
we need only show that $\partial y^i/\partial x^j\in C^{0,\alpha}$.

Differentiating (\ref{ode}) with respect to $x^j$ yields
\begin{displaymath}
\frac {\partial^2 y^i}{\partial x^j\partial t}(x,t) =
\frac {\partial V^i}{\partial x^k} (y(x,t))
\frac {\partial y^k}{\partial x^j}(x,t).
\end{displaymath}
Write $A^i_k(x,t) = \partial V^i/\partial x^k(y(x,t))$ and $u^i_j(x,t)
=\partial y^i/\partial x^j(x,t)$.  Since the composition of a $C^1$
function and a $C^{0,\alpha}$ function is in $C^{0,\alpha}$, $A^i_k$ is of
class $C^{0,\alpha}$.  Since $\partial y^i/\partial x^j\partial t$ is
continuous, we can commute the second derivatives and conclude that $u^i_j$
satisfies the {\it linear} initial value problem
\begin{align*}
\frac{\partial u^i_j}{\partial t}(x,t) &= A^i_k(x,t) u^k_j(x,t),\\
u^i_j(x,0) &= \delta^i_j.
\end{align*}
(The initial condition follows from differentiating (\ref{iv}) with respect
to $x^j$.)

For any fixed $x_1$ and $x_2$, define
\begin{displaymath}
v^i_j(t) = u^i_j(x_1,t) - u^i_j(x_2,t).
\end{displaymath}
Note that $v^i_j(0) = 0$, and that the equation for $u^i_j$ implies
\begin{displaymath}
\left| \frac {\partial v}{\partial t}(t)\right|
\le \|u\|_{L^\infty}  \| A \| _{C^{0,\alpha}} |x_1-x_2|^\alpha +
\|A\| _{L^\infty} |v(t)|.
\end{displaymath}
By one version of a standard ODE comparison theorem \cite[Cor.\ III.4.3,
p.\ 27]{Hartman}, this implies that $|v(t)|$ is bounded by the solution to
the IVP
\begin{align*}
y'(t) &=   \| A \| _{C^{0,\alpha}}\left( \|u\|_{L^\infty} |x_1-x_2|^\alpha  +
y(t) \right),\\
y(0)&=0,
\end{align*}
which is
\begin{displaymath}
y(t) = (\exp t\|A\|_{C^{0,\alpha}} - 1) \|u\|_{L^\infty} |x_1-x_2|^\alpha.
\end{displaymath}
This yields the desired estimate.
\end{pf}

\noindent{\bf Remark.}
This proof can easily be adapted to obtain higher H\"older regularity for
the solution if $V\in C^{m,\alpha}$, $m > 1$.  However, as was pointed out
to me by Curt McMullen, the analogous result fails for vector fields that
are merely $C^{0,\alpha}$; in fact, there are vector fields with
coefficients in the ``Zygmund class'', thus in $C^{0,\alpha}$ for
every $0<\alpha<1$ and uniquely integrable, but whose solutions have
arbitrarily bad H\"older exponents.  See \cite{Reimann} for further
discussion.

\bigskip
\noindent\fbox{\sc Version 1.0 -- September 19, 1994}

\end{document}